\documentclass[12pt]{amsart}

\newtheorem{THM}{Theorem}
\newtheorem{LMA}[THM]{Lemma}

\usepackage{amsmath}
\usepackage{amssymb}
\usepackage{euscript}
\usepackage[pdftex]{graphicx}
\newcommand{\pix}[3]{\begin{figure}[tbp]
\centering{\includegraphics[width=#1\textwidth]{#2}} 
\caption{#3}
\end{figure}}

\newcommand{\showon}{\begin{eqnarray*}}
\newcommand{\showoff}{\end{eqnarray*}}
\newcommand{\Showon}{\begin{eqnarray}}
\newcommand{\Showoff}{\end{eqnarray}}

\newcommand{\none}{\varnothing}
\newcommand{\drop}{\smallsetminus}

\newcommand{\goesto}{\rightarrow}
\newcommand{\one}{\boldsymbol{1}}
\newcommand{\zero}{\boldsymbol{0}}
\newcommand{\LHS}{\mathrm{LHS}}
\newcommand{\RHS}{\mathrm{RHS}}
\newcommand{\LDIFF}{\mathrm{L}_0}
\newcommand{\RDIFF}{\mathrm{R}_0}

\newcommand{\T}{\EuScript{T}}

\newcommand{\X}{\EuScript{X}} 

\newcommand{\Y}{\EuScript{Y}} \newcommand{\y}{\mathbf{y}}

\begin{document}

\title[Rayleigh monotonicity]{A combinatorial proof of Rayleigh monotonicity for graphs}

\author{J. Cibulka}
\author{J. Hladk\'y}
\thanks{Research of JC and JH supported by the DIMACS-DIMATIA REU grant (NSF CNS
0138973)}
\address{Department of Applied Mathematics\\
Charles University\\
Malostransk\'e n\'am. 25\\
118 00 Praha 1\\
Czech Republic}
\email{\texttt{cibulka@kam.mff.cuni.cz}}
\email{\texttt{hladky@kam.mff.cuni.cz}}
\author{M.A. LaCroix}
\author{D.G. Wagner}
\thanks{Research of DGW supported by NSERC Discovery Grant OGP0105392.}
\address{Department of Combinatorics and Optimization\\
University of Waterloo\\
Waterloo, Ontario, Canada\ \ N2L 3G1}
\email{\texttt{malacroix@math.uwaterloo.ca}}
\email{\texttt{dgwagner@math.uwaterloo.ca}}


\begin{abstract}
We give an elementary, self-contained, and purely combinatorial proof of
the Rayleigh monotonicity property of graphs.
\end{abstract}

\maketitle

Consider a (linear, resistive) electrical network -- this is a connected
graph $G=(V,E)$ and a set of positive real numbers $\y=\{y_e:\ e\in E\}$ indexed by $E$.
In this paper we allow graphs to have loops and/or multiple edges.
The value of $y_e$ is interpreted as the
electrical conductance of a wire joining the vertices incident with $e$.
For any edge $e\in E$, there is a simple formula for the effective conductance
$\Y_e(G;\y)$ of the rest of the graph $G\drop \{e\}$, measured between the ends of $e$.
This is due to Kirchhoff \cite{Ki} and is also known as Maxwell's Rule \cite{Max}.
For a subset $S\subseteq E$, let
$$\y^S = \prod_{c\in S} y_c.$$
Spanning subgraphs of $G$ will be identified naturally with their edge-sets.
Let $\T(G)$ be the set of all spanning trees of $G$, and let
$$T(G;\y) = \sum_{T\in\T(G)} \y^T$$
be the \emph{tree-generating polynomial} of $G$.  Kirchhoff's Formula
for the effective conductance of an electrical network is
$$\Y_e(G;\y)=\frac{T(G\drop \{e\};\y)}{T(G/\{e\};\y)}.$$

It is physically intuitive that if all edge-conductances are positive and
we increase the conductance of an edge $f$, then the effective conductance
$\Y_e(G;\y)$ does not decrease.  That is,
$$\frac{\partial}{\partial y_f} \frac{T(G\drop \{e\};\y)}{T(G/\{e\};\y)} \geq 0$$
provided that $y_c>0$ for all $c\in E$.  This property is known as
\emph{Rayleigh monotonicity}.  To keep formulas readable, we often suppress
the variables $\y$ (and sometimes the graph $G$)
from the notation unless they require particular attention.
A further shorthand is to write
$$T(G) = T^e(G) + y_e T_e(G),$$
in which $T^e(G) = T(G\drop \{e\})$ and $T_e(G) = T(G/\{e\})$.
Applying the quotient rule and a little algebra, Rayleigh monotonicity
is seen to be equivalent to the inequality
$$T_e^f(G)T_f^e(G) - T_{ef}(G)T^{ef}(G) \geq 0$$
whenever $y_c>0$ for all $c\in E$.  We also use the notation
$\T_e^f(G)=\T((G\drop\{f\})/\{e\})$ for the set of spanning trees of $(G\drop\{f\})/\{e\}$,
\emph{et cetera}.

In fact, the Rayleigh monotonicity property of graphs follows from a
more precise -- and rather surprising -- combinatorial identity.
Fix two distinct edges $e,f\in E$, and orient them arbitrarily.
Let $\X=\X(G;e,f)$ denote the set of spanning forests $F$ of $G$
such that both $F\cup\{e\}$ and $F\cup\{f\}$ are spanning trees of $G$.
Thus, $F\cup\{e,f\}$ contains a unique cycle $C$, and $C$ contains
both $e$ and $f$.  Let $\X^{+}=\X^+(G;e,f)$ denote the set of those $F\in\X$ for
which the edges $e$ and $f$ are oriented in the same direction around
the corresponding cycle $C$.  Let $\X^{-}=\X^-(G;e,f)$ denote the set of those
$F\in\X$ for which the edges $e$ and $f$ are oriented in opposite
directions around the corresponding cycle $C$.  Define
$$X^{+}(G;e,f) = \sum_{F\in\X^{+}(G;e,f)} \y^{F}$$
and
$$X^{-}(G;e,f) = \sum_{F\in\X^{-}(G;e,f)} \y^{F}.$$

\begin{THM}
Let $G=(V,E)$ be a connected graph, and let $\y=\{y_c:\ c\in E\}$
be indeterminates.
With the notation given above, for distinct $e,f\in E$,
$$T_{e}^{f}(G) T_{f}^{e}(G) - T_{ef}(G) T^{ef}(G) =
\left[ X^{+}(G;e,f) - X^{-}(G;e,f) \right]^{2}.$$
\end{THM}
The case of Theorem 1 with all $\y\equiv\one$ appears as equation (2.34) of Brooks, Smith, Stone,
and Tutte \cite{BSST}, and the generalization of this to regular matroids is Theorem 2.1 of Feder
and Mihail \cite{FM}.  The case of general conductances $\y>\zero$ on graphs can be found in Section 3.8
of Balabanian and Bickart \cite{BB}.  Choe \cite{Ch1,Ch2,Ch3} gives two proofs of Theorem~1.
One -- based on Tutte's theory of chain groups as in \cite{FM} -- generalizes to all
sixth-root-of-unity matroids but gives a less precise description of the right-hand side.
The other proof (as in \cite{BB,BSST}) uses the All-Minors Matrix-Tree Theorem \cite{Cha}
and Jacobi's theorem on complementary minors of inverse matrices,
together with substantial and elaborate algebraic manipulations.  Neither of these
proofs is completely combinatorial.  Here we give a proof of Theorem 1
that is elementary, self-contained, and purely combinatorial.  This is not a
direct bijective proof, however:\  we proceed by induction on the number of edges
of the graph, and in the induction step we resort to natural 2:2 or 2:1 correspondences
(as well as 1:1 correspondences) and in one case we employ a sign-reversing involution.
Similar versions of this proof were found independently -- by JC and JH, and by MAL and
DGW -- at about the same time.  An earlier description is given by Cibulka and Hladk\'y \cite{CH}.
We thank Bill Jackson for his careful reading of an earlier version of this argument that suffered from
an inaccurate description of the last case.\\

\noindent\textbf{\emph{Proof of Theorem 1.}}

Consider a connected graph $G=(V,E)$ and conductances $\y=\{y_c:\ c\in E\}$ and a pair
of edges $e,f\in E$.  For short, let us write $T_e^f$ for $T_e^f(G;\y)$ and
$X^+$ for $X^+(G;e,f;\y)$, and so on.  We establish Theorem 1 in the form
\begin{eqnarray}
T_e^f T_f^e + 2 (X^+) (X^-) &=& T_{ef} T^{ef} + (X^+)^2 + (X^-)^2
\end{eqnarray}
by induction on the number of edges of $G$. 
To establish the polynomial equation (1), we consider an arbitrary monomial $\y^\alpha$ and show
that the coefficients of $\y^\alpha$ on each side of the equation are equal.
A general term appearing in equation (1) is indexed by a pair of sets
$(A,B)$ with $A,B\subseteq E$, and it contributes the monomial $\y^A \y^B$.
The pairs in the set
$$(\T_e^f\times\T_f^e) \cup (\X^+\times\X^-) \cup (\X^-\times \X^+)$$
contribute to the left-hand side of (1), and the pairs in the set
$$(\T_{ef}\times\T^{ef}) \cup (\X^+\times\X^+) \cup (\X^-\times \X^-)$$
contribute to the right-hand side of (1).  The \emph{type} of a pair $(A,B)$
that contributes to equation (1) is defined to be that one of the six expressions
$$
\T_e^f\T_f^e,\ \X^+\X^-,\ \X^-\X^+,\ \T_{ef}\T^{ef},\ \X^+\X^+,\ \X^-\X^-
$$
which describes the set which contains it.

Notice that $T$ and $X^+$ and $X^-$ are \emph{multiaffine} -- each indeterminate
$y_c$ occurs at most to the first power.  Therefore, equation (1) is
at most quadratic in each variable.  Thus, we need only consider 
monomials $\y^\alpha$ such that $\alpha:E\goesto\{0,1,2\}$.  Notice also
that neither $y_e$ nor $y_f$ occur in equation (1), so we need only consider
monomials $\y^\alpha$ such that $\alpha(e)=\alpha(f)=0$.  Furthermore, if
$g\in E$ is a loop in $G$ then $y_g$ does not occur in equation (1).  Thus,
it suffices to prove Theorem 1 for loopless graphs $G$.  We even have the following
more substantial reduction.

\begin{LMA}
Let $G$ be a connected graph that is the union of subgraphs $H$ and $J$ which have
exactly one vertex (and no edges) in common.  If (1) holds for $H$ and for $J$, for
every choice of edges $e,f$, then (1) holds for $G$, for every choice of edges $e,f$.
\end{LMA}
\begin{proof}
The key observation is that $T(G)=T(H)T(J)$.
Consider any two distinct edges $e,f$ of $G$.  Up to symmetry, there are two cases:\
either $e,f$ are both in $H$, or $e$ is in $H$ and $f$ is in $J$.

If $e,f$ are both in $H$, then $X^+(G)=X^+(H)T(J)$ and $X^-(G)=X^-(H)T(J)$, so
that equation (1) for $(G,e,f)$ is just $T(J)^2$ times equation (1) for $(H,e,f)$.
By the hypothesis, this equation holds.

If $e$ is in $H$ and $f$ is in $J$, then $X^+(G)=X^-(G)=0$ and
$T_e^f(G)=T_e(H)T^f(J)$, $T_f^e(G)=T^e(H)T_f(J)$,
$T_{ef}(G)=T_e(H)T_f(J)$, and $T^{ef}(G)=T^e(H)T^f(J)$.
In this case, equation (1) for $(G,e,f)$ can be verified directly.
\end{proof}

\noindent\textbf{Basis of Induction.}

The basis of induction consists of three cases.

\underline{Case 1:\ $G=C_2$.}
The graph consists of two edges $e,f$ in parallel.  Orient them to
form a directed $2$-cycle.  We have
$T_e^f=T_f^e=1$, $T_{ef}=T^{ef}=0$, $X^+=1$ and $X^-=0$.  Equation (1)
states that
$$ 1\cdot 1 + 1\cdot 0 + 0\cdot 1 =
0\cdot 0 + 1\cdot 1 + 0\cdot 0, $$
which is clearly true.

\underline{Case 2:\ $G=P_3$.}
The graph consists of two edges $e,f$ incident at one common vertex.
Direct them arbitrarily.  We have
$$T_e^f=T_f^e=T^{ef}=X^+=X^-=0\ \ \mathrm{and}\ \ T_{ef}=1$$
so that both sides of equation (1) are zero.

\underline{Case 3:\ $G=K_4$.}
As will be seen, our induction step does not apply to $K_4$.
(The induction step in \cite{CH} does, however.)
Thus, we verify equation (1) for $K_4$ directly.
Consider $K_4$ with edges labelled as in Figure 1.
\pix{1.0}{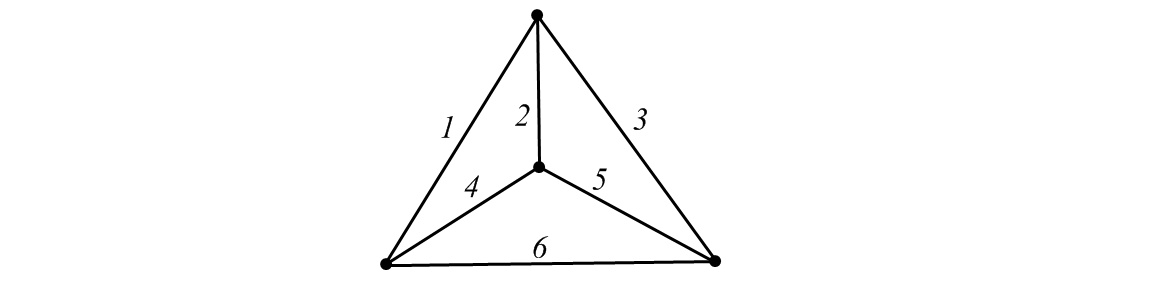}{The basis of induction $K_4$.}
The edges $e$ and $f$ are either adjacent or not, so
(up to automorphism) we have the following two subcases.\\
$\bullet$\ \ $\{e,f\}=\{1,2\}$.  Direct edges $1$ and $2$ towards their common vertex.  Then
\showon
T_1^2  &=& y_3 y_4 + y_3 y_5 + y_4 y_5 + y_4 y_6 + y_5 y_6 \\ 
T_2^1  &=& y_3 y_4 + y_3 y_6 + y_4 y_5 + y_4 y_6 + y_5 y_6 \\
T_{12} &=& y_3 + y_5 + y_6 \\
T^{12} &=& y_3(y_4 y_5 + y_4 y_6+ y_5 y_6) \\
X^+    &=& 0 \\
X^-    &=& y_3 y_4 + y_4 y_5 + y_4 y_6 + y_5 y_6
\showoff
$\bullet$\ \ $\{e,f\}=\{1,5\}$. Direct edges $1$ and $5$ towards the vertices they share
with the edge $2$.  Then
\showon
T_1^5  &=&  (y_2 + y_4)(y_3 + y_6)\\
T_5^1  &=&  (y_2 + y_3)(y_4 + y_6)\\
T_{15} &=&  y_2 + y_3 + y_4 + y_6\\
T^{15} &=&  y_2y_3y_4 + y_2y_3y_6 + y_2y_4y_6 + y_3y_4y_6\\
X^+    &=&  y_3y_4\\
X^-    &=&  y_2y_6
\showoff
In each subcase equation (1) is easily verified.\\

\noindent\textbf{Induction Step.}

Consider a loopless connected graph $G=(V,E)$ with $n$ vertices and $m\geq 3$
edges, let $\y=\{y_c:\ e\in E\}$, and fix distinct edges $e,f\in E$.  Assume
that equation (1) holds for any connected graph with at most $m-1$ edges.
By Lemma 2 we can assume that $G$ has no cut-vertices, and hence no cut-edges.

Consider any monomial $\y^\alpha$ such that $\alpha:E\goesto\{0,1,2\}$ and
$\alpha(e)=\alpha(f)=0$.  We need to show that
\begin{eqnarray}
[\y^\alpha](T_e^f T_f^e + 2 (X^+) (X^-) )
=
[\y^\alpha](T_{ef} T^{ef} + (X^+)^2 + (X^-)^2)
\end{eqnarray}
in which $[\y^\alpha]P(\y)$ denotes the coefficient of the monomial $\y^\alpha$
in the polynomial $P(\y)$.

Now, two easy reductions.  If $\alpha(g)=0$ for some $g\in E\drop\{e,f\}$,
then we can use the fact that a pair
$(A,B)$ contributes $\y^A \y^B = \y^\alpha$ to one side of equation (1)
for $G$ if and only if it contributes $\y^A \y^B = \y^\alpha$ to the same
side of equation (1) for $G\drop \{g\}$.  By induction we can assume that (1)
holds for $G\drop \{g\}$ (since $g$ is not a cut-edge), and we conclude that (2) holds
for $G$ and $\alpha$.

Similarly, if $\alpha(g)=2$ for some $g\in E\drop \{e,f\}$, then
we can use the fact that a pair $(A,B)$ contributes $\y^A \y^B =
\y^\alpha$ to one side of equation (1) for $G$ if and only if
$(A\drop \{g\}, B\drop \{g\})$ contributes $y_g^{-2}\y^\alpha$ to the same
side of equation (1) for $G/g$.  By induction we can assume that (1)
holds for $G/g$, and we conclude that (2) holds for $G$ and $\alpha$.

For the rest of the induction step we need only consider the
monomial $\y^\gamma$ such that $\gamma(e)=\gamma(f)=0$ and
$\gamma(g)=1$ for all $g\in E\drop\{e,f\}$.
It will be convenient to have the notation
\showon
\LHS(G) &=& [\y^\gamma](T_e^f T_f^e + 2 (X^+) (X^-) )\\
\mathrm{and}\ \ \ \
\RHS(G) &=& [\y^\gamma](T_{ef} T^{ef} + (X^+)^2 + (X^-)^2)
\showoff
(since the edges $e$ and $f$ remain fixed throughout, we suppress them
from the notation).  To complete the induction step it suffices to prove that
\begin{eqnarray}
\LHS(G) = \RHS(G).
\end{eqnarray}
The pairs $(A,B)$ contributing to equation (3) are ordered partitions of
the edge-set $E\drop\{e,f\}$ into two (possibly empty) subsets.

\begin{LMA}
Let  $G=(V,E)$  be a graph, let  $e,f\in E$,  and let  $K$  be an edge-cut of
$G$  disjoint from  $\{e,f\}$.  For any pair  $(A,B)$  contributing to
equation (3),  both  $A\cap K\neq\none$  and $B\cap K\neq\none$.
\end{LMA}
\begin{proof}
For any pair $(A,B)$ contributing to equation (3), the
spanning subgraphs $A\cup\{e,f\}$ and $B\cup\{e,f\}$ of $G$
are both connected.  Thus, if $K\cap\{e,f\}=\none$ then both
$A\cap K\neq\none$ and $B\cap K\neq\none$.
\end{proof}

The polynomial equation (1) is homogeneous of degree $2(n-2)$.
The monomial $\y^\gamma$ has degree $m-2$.  Thus, equation (3) is
trivial except in the case that $m=2n-2$, and so we reduce to this case.
Since the sum of the degrees of the vertices of $G$ is $2m=4n-4$,
it follows that $G$ must contain a vertex of degree at most three.
The following cases are indexed by the degree of a vertex at which
we perform a reduction of the graph $G$, in order to apply the induction
hypothesis.\\

\noindent\textbf{Cases 0 and 1.}

A vertex of degree zero in $G$ is not possible since in the induction
hypothesis we assume that the graph is connected and loopless with
$m\geq 3$ edges.  If $G$ had a vertex of degree one then, since $m\geq 3$
and $G$ is loopless, $G$ would contain a cut-vertex.  Since we have reduced
to the case that $G$ has no cut-vertices, a vertex of degree one in $G$ is
also impossible.\\ 

\noindent\textbf{Case 2.}

If $v_*$ is a vertex of degree $2$ in $G$, then there are three subcases:\\
(i)  \ \ $v_*$ is incident with neither $e$ nor $f$;\\
(ii) \ \ $v_*$ is incident with $e$ but not $f$;\\
(iii)\ \ $v_*$ is incident with both $e$ and $f$.\\
(By symmetry, (ii) also covers the case that $v_*$ is incident with $f$ but not $e$.)

$\bullet$\ \
In subcase (i) let $v_*$ be incident with $g$ and $h$.
Then
$$\LHS(G)=2\cdot\LHS(G\drop \{v_*\})$$
since each pair $(A,B)$ contributing to $\LHS(G\drop \{v_*\})$
gives rise to two pairs $(A\cup \{g\},B\cup \{h\})$ and $(A\cup \{h\},B\cup \{g\})$
contributing to $\LHS(G)$.  (By applying Lemma 3 with $K=\{g,h\}$ one
sees that these are the only possibilities.)
Similarly,
$$\RHS(G)=2\cdot\RHS(G\drop \{v_*\}).$$
Note that $G\drop \{v_*\}$ is connected, since $G$ has no cut-vertices.
By the induction hypothesis applied to $G\drop \{v_*\}$, this suffices
to establish equation (3) in this subcase.

$\bullet$\ \
In subcase (ii) let $v_*$ be incident with $e$ and $g$.
Every forest $F$ contributing $\y^F$ to $X^+(G)$ or $X^-(G)$ must use the
edge $g$.  Therefore,
$$[\y^\gamma](X^+)(X^-)=[\y^\gamma](X^+)^2=[\y^\gamma](X^-)^2=0.$$
Also, if $T\in \T^e(G)$ then $g\in T$.  Therefore,
\showon
[\y^\gamma]T_e^f(G)T_f^e(G)
&=&
[\y^\gamma] y_g T_e^{fg}(G)T_{fg}^e(G)\\
&=&
[\y^\gamma] y_g T_{ef}^{g}(G)T_g^{ef}(G)
=[\y^\gamma]T_{ef}(G)T^{ef}(G),
\showoff
in which the second equality is a consequence of the 1:1 correspondence
$(A,B)\leftrightarrow (B\cup \{e\}\drop \{g\},A\cup \{g\}\drop \{e\})$.  This proves
equation (3) directly in this subcase.

$\bullet$\ \
In subcase (iii) vertex $v_*$ is incident with $e$ and $f$.
Re-orienting $e$ and/or $f$, if necessary, we may assume that
in the definition of $X^+$ and $X^-$ both $e$ and $f$ are
directed towards $v_*$.  Notice that $T_e^f(G)=T_f^e(G)=T(G\drop \{v_*\})$
and $T^{ef}(G)=0$, and $X^+(G)=0$ and $X^-(G)=T(G\drop \{v_*\})$.  Thus,
in this subcase, equation (3) reduces to 
$$
\LHS(G) = [\y^\gamma]T(G\drop \{v_*\})^2 = \RHS(G),
$$
completing the analysis of Case 2.\\

For the rest of the induction step we need only consider
a two-connected graph $G$ with $n$ vertices
and $m=2n-2\geq 3$ edges, and with minimum degree three.
Such a graph must have at least four vertices of degree three.
In fact, one of the two following cases must occur:\\
\underline{Case 3:}\ \ there is a $3$-valent vertex $v_*$ incident with neither $e$ nor $f$;\\
\underline{Case 4:}\ \ $e$ and $f$ are not adjacent, their four ends are $3$-valent,
and every other vertex of $G$ is $4$-valent.\\

\noindent\textbf{Case 3.}

Let $v_{*}$ be a vertex of degree three in
$G$ that is not incident with either edge $e$ or edge $f$.  Let
$a=\{v_*,a_*\}$,
$b=\{v_*,b_*\}$, and
$c=\{v_*,c_*\}$  be the edges incident with $v_{*}$ in $G$.
Consider three new edges (not in $G$) with ends given by
$ab = \{a_*,b_*\}$ and
$bc = \{b_*,c_*\}$ and
$ac = \{a_*,c_*\}$.
Form the graphs
$H(ab)=(G\drop \{v_{*}\})\cup\{ab\}$,
$H(bc)=(G\drop \{v_{*}\})\cup\{bc\}$, and
$H(ac)=(G\drop \{v_{*}\})\cup\{ac\}$, as
depicted in Figure 2.
(Note that $a_*=b_*$ is possible, for example, in which case $ab$ is a loop.)
Since $v_*$ is not a cut-vertex of $G$, it follows that all of the graphs $H(ab)$,
$H(ac)$, and $H(bc)$ are connected.
\pix{1.0}{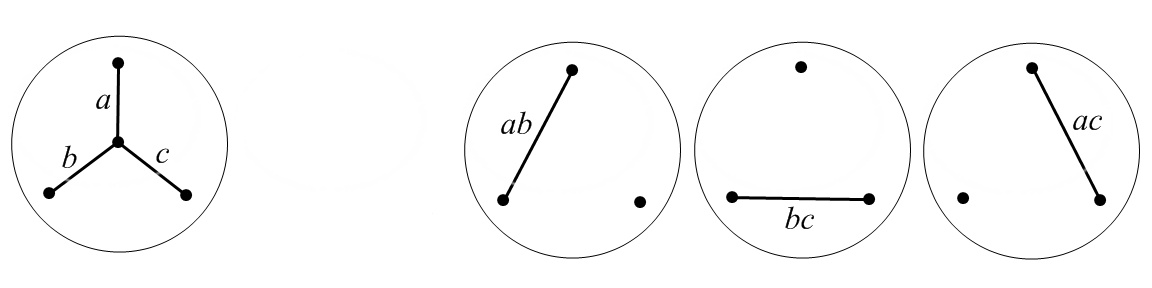}{Reduction at a $3$-valent vertex.}

\begin{table}
\begin{center}
\begin{tabular}{|c|c|c|c|}\hline
&&& \\
$(A\cap K, B\cap K)$ & $H$ & $A'$ & $B'$ \\ \hline
$(\{a,b\},\{c\})$ & $H(ab)$ & $(A\drop \{v_*\})\cup\{ab\}$ & $B\drop \{v_*\}$ \\ 
$(\{c\},\{a,b\})$ & $H(ab)$ & $A\drop \{v_*\}$ & $(B\drop \{v_*\})\cup\{ab\}$ \\ \hline
$(\{a,c\},\{b\})$ & $H(ac)$ & $(A\drop \{v_*\})\cup\{ac\}$ & $B\drop \{v_*\}$ \\
$(\{b\},\{a,c\})$ & $H(ac)$ & $A\drop \{v_*\}$ & $(B\drop \{v_*\})\cup\{ac\}$ \\ \hline
$(\{b,c\},\{a\})$ & $H(bc)$ & $(A\drop \{v_*\})\cup\{bc\}$ & $B\drop \{v_*\}$ \\
$(\{a\},\{b,c\})$ & $H(bc)$ & $A\drop \{v_*\}$ & $(B\drop \{v_*\})\cup\{bc\}$ \\ \hline
\end{tabular}
\end{center}
\caption{1:1 in Case 3.}
\end{table}

We claim that
\Showon
\begin{array}{l}
\LHS(G) = \\
\LHS(H(ab))+\LHS(H(ac))+\LHS(H(bc))
\end{array}
\Showoff
and
\Showon
\begin{array}{l}
\RHS(G) = \\
\RHS(H(ab))+\RHS(H(ac))+\RHS(H(bc)).
\end{array}
\Showoff
By the induction hypothesis, equation (3) holds for each of
the graphs $H(ab)$, $H(ac)$, and $H(bc)$, so that (4) and (5)
suffice to prove equation (3) for $G$.

Consider any pair $(A,B)$ contributing to equation (3) for $G$.
By Lemma 3 applied with $K=\{a,b,c\}$, the induced ordered partition
$(A\cap K, B\cap K)$ of $\{a,b,c\}$ is one of the six cases presented in
the first column of Table 1.  The corresponding entry of the second
column indicates the graph $H$ to which it is assigned.  The corresponding
entries of the third and fourth columns describe the pair $(A',B')$ associated 
with $(A,B)$ that contributes to equation (3) for $H$.  It is easy to see that
the pairs $(A,B)$ and $(A',B')$ have the same type.
Thus, the construction described in Table 1 gives a 1:1 correspondence
$(A,B)\leftrightarrow(A',B')$ that proves equations (4) and (5).
Therefore, we have verified equation (3) for $G$ in Case 3.\\

\noindent\textbf{Case 4.}

If $G=K_4$ then we have already verified the conclusion as part of the basis of
induction.  Otherwise, let $v_{*}$ be a vertex of
$G$ that is not incident with either edge $e$ or edge $f$, and note that
$v_*$ has degree $4$ and is not a cut-vertex.  Let
$a=\{v_*,a_*\}$,
$b=\{v_*,b_*\}$, 
$c=\{v_*,c_*\}$, and
$d=\{v_*,d_*\}$  be the edges incident with $v_{*}$ in $G$.
Consider six new edges (not in $G$) with ends given by
$ab = \{a_*,b_*\}$,
$ac = \{a_*,c_*\}$,
$ad = \{a_*,d_*\}$,
$bc = \{b_*,c_*\}$,
$bd = \{b_*,d_*\}$, and
$cd = \{c_*,d_*\}$.
Form the graphs
$H(ab|cd)=(G\drop \{v_{*}\})\cup\{ab,cd\}$,
$H(ac|bd)=(G\drop \{v_{*}\})\cup\{ac,bd\}$, and
$H(ad|bc)=(G\drop \{v_{*}\})\cup\{ad,bc\}$, as
depicted in Figure 3.
(Note that $a_*=b_*$ is possible, for example, in which case $ab$ is a loop.)
Since $v_*$ is not a cut-vertex, all three of
$H(ab|cd)$, $H(ac|bd)$, and $H(ad|bc)$ are connected.
By induction we can assume that equation (3) holds for each of these graphs,
which we will call \emph{$H$ graphs}.
\pix{1.0}{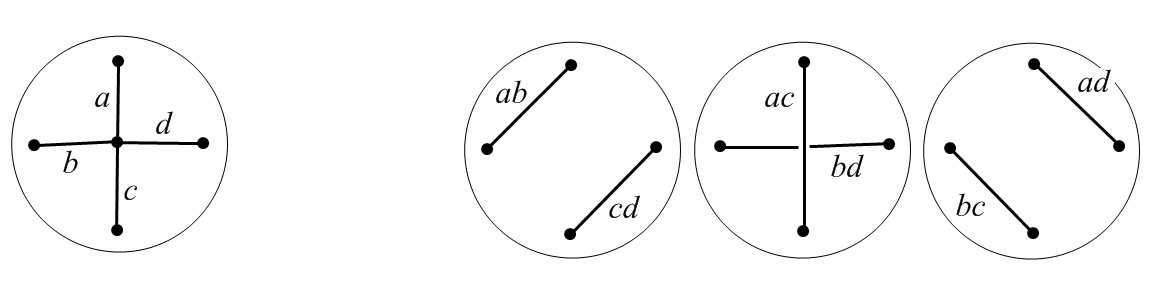}{Reduction at a $4$-valent vertex.}

The analogues of (4) and (5) do not hold in this case.
Instead, we define numbers $\LDIFF$ and $\RDIFF$ by
\Showon
\begin{array}{l}
\LHS(G) + \LDIFF = \\
\LHS(H(ab|cd))+\LHS(H(ac|bd))+\LHS(H(ad|bc))
\end{array}
\Showoff
and
\Showon
\begin{array}{l}
\RHS(G) + \RDIFF = \\
\RHS(H(ab|cd))+\RHS(H(ac|bd))+\RHS(H(ad|bc))
\end{array}
\Showoff
and complete the proof by showing that
\Showon
\LDIFF=\RDIFF.
\Showoff
The induction hypotheses (3) for the three $H$ graphs, together with equations (6), (7),
and (8), suffice to prove equation (3) for $G$.

To prove equation (8) we compare the set of pairs $(A,B)$ contributing to equation (3)
for $G$ with the set of pairs $(A',B')$ contributing to equation (3) for 
the three $H$ graphs.  This results in a combinatorial description of
$\LDIFF$ and $\RDIFF$, from which equation (8) follows easily.  To make
this comparison, consider any pair $(A,B)$ contributing to equation (3) for $G$. 
By Lemma 3 applied with $K=\{a,b,c,d\}$, the induced unordered partition
$\{A\cap K, B\cap K\}$ of $\{a,b,c,d\}$ falls into one of the following two subcases: either\\
(i)\ \   $\{A\cap K, B\cap K\}$ has two blocks of size two, or\\
(ii)\ \ $\{A\cap K, B\cap K\}$  has one block of size three and one block of size one.

\underline{In subcase (i)},\ if $(A\cap K,B\cap K) = (\{h,i\},\{j,k\})$ for some labelling
$\{h,i,j,k\}=\{a,b,c,d\}$, then let
$$(A',B')=((A\drop \{v_*\})\cup \{hi\}, (B\drop \{v_*\})\cup \{jk\})$$
as depicted in Figure 4.
\pix{1.0}{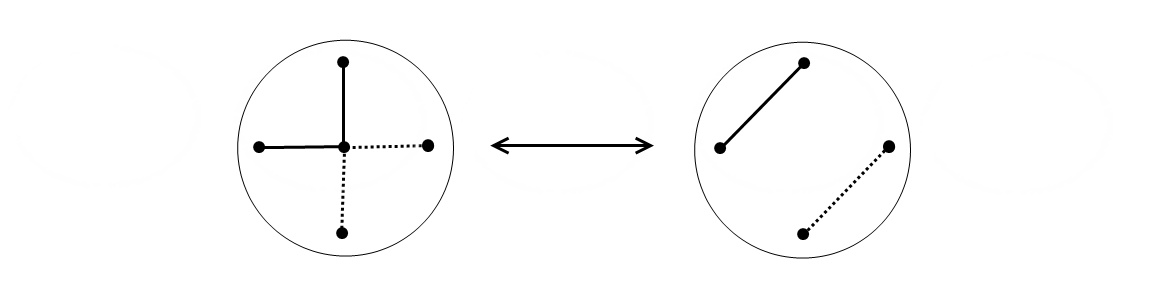}{1:1 in Case 4(i).}
This gives a 1:1 correspondence $(A,B)\leftrightarrow(A',B')$
between the pairs $(A,B)$ for $G$ in subcase (i) and the pairs $(A',B')$ for some
$H(hi|jk)$ with edges $hi$ and $jk$ in different sets from the pair $(A',B')$.
It is clear that the pairs $(A,B)$ and $(A',B')$ are of the same type.

\underline{In subcase (ii)},\ we describe in detail the situation in which $A\cap K=\{h,i,j\}$ and
$B\cap K=\{k\}$ for some labelling $\{h,i,j,k\}=\{a,b,c,d\}$.
(The other situation is analogous under exchange of $A$ with $B$.)
Note that the vertices $h_*$, $i_*$, and $j_*$ are distinct, since $A$ is a forest.
We distinguish among a number of possibilities, exactly one of which occurs:
either\\
(ii-a)\ \ both $A$ and $B$ are spanning trees, or\\
(ii-b)\ \ both $A$ and $B$ are forests with two components.\\
In case (ii-b), both $A\cup \{e\}$ and $A\cup \{f\}$ are spanning trees.  Let $A[e]$ be the unique path
in $A\cup \{e\}$ from $v_*$ to $k_*$, and let $A[f]$ be the unique path in $A\cup \{f\}$ from $v_*$ to $k_*$.
Let $\widetilde{e}$ be the edge incident with $v_*$ in $A[e]$, and let $\widetilde{f}$ be the edge
incident with $v_*$ in $A[f]$.  We further divide case (ii-b) as follows: either\\
(ii-b1)\ \ the edges $\widetilde{e}$ and $\widetilde{f}$ are equal, or\\
(ii-b2)\ \ the edges $\widetilde{e}$ and $\widetilde{f}$ are not equal.

It remains to associate to each such pair $(A,B)$, a pair
$(A',B')$ that contributes to equation (3) for one of the $H$ graphs.
In the situation we are describing in detail
($A\cap K=\{h,i,j\}$ and $B\cap K=\{k\}$) we always put $B'=B\drop \{v_*\}$.

In case (ii-a), let $h_*$ be the vertex
adjacent to $v_*$ on the unique path from $v_*$ to $k_*$ in $A$, and put
$\widetilde{e}=\widetilde{f}=h$.  We associate two different pairs $(A',B')$
and $(A'',B')$ with $(A,B)$ in this case, by putting
$$A'=(A\drop \{v_*\})\cup\{hi,jk\}
\ \ \ \mathrm{and}\ \ \
A''=(A\drop \{v_*\})\cup\{hj,ik\}.$$
Figure 5 illustrates this construction (with $k=d$ and $h=b$).
Note that the pairs $(A,B)$, $(A',B')$, and $(A'',B')$ have the same type.
\pix{1.0}{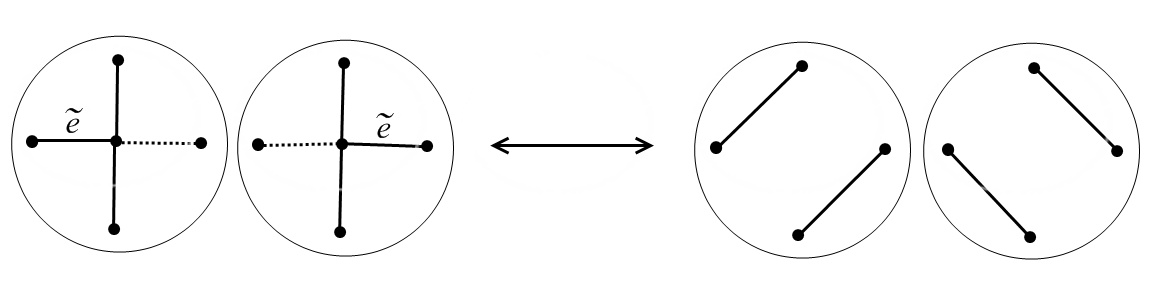}{2:2 in Cases 4(ii-a) and 4(ii-b1).}
Consider the pair of sets $(P,Q)$ with $P=A\cup \{k\}\drop \{h\}$ and $Q=B\cup \{h\}\drop \{k\}$.
This pair contributes to equation (3) for $G$, and is also in case (ii-a).
Moreover, the two pairs $(P',Q')$ and $(P'',Q')$ associated to it are
in fact the same as $(A',B')$ and $(A'',B')$.

In case (ii-b1), let $h_*$ be the neighbour of $v_*$ that is 
incident with $\widetilde{e}=\widetilde{f}$.  In this case we apply the same construction
as in case (ii-a).
Note that the pairs $(A,B)$, $(A',B')$, and $(A'',B')$ have the same type.
The pair of sets $(P,Q)$ with $P=A\cup \{k\}\drop \{h\}$ and $Q=B\cup \{h\}\drop \{k\}$
contributes to equation (3) for $G$, and is also in case (ii-b1), and
the two pairs $(P',Q')$ and $(P'',Q')$ associated to it are
the same as $(A',B')$ and $(A'',B')$.

In case (ii-b2) let $h_*$ be the neighbour of $v_*$ that is
incident with $\widetilde{e}$, and let $i_*$ be the neighbour of $v_*$ that
is incident with $\widetilde{f}$.
In this case we associate a single pair $(A',B')$ with $(A,B)$ by putting
$$A'=(A\drop \{v_*\})\cup\{hi,jk\}$$
as depicted in Figure 6 (with $k=d$, $h=a$ and $i=c$).
Note that the pairs $(A,B)$ and $(A',B')$ have the same type.
\pix{1.0}{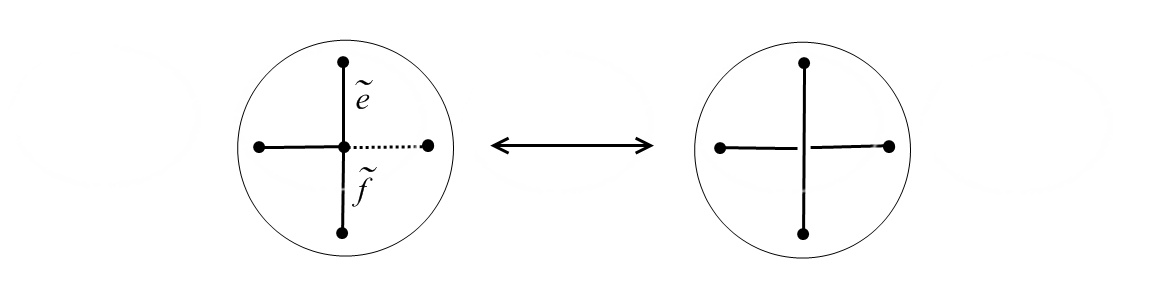}{1:1 in Case 4(ii-b2).}

Table 2 summarizes these correspondences in Case 4(ii).
\begin{table}
\begin{center}
\begin{tabular}{|c|c|c|c|}\hline
&&& \\
$G$ case & $m$:$m$ & $H$ case & $A'$ \\ \hline
(a)  & 2:2 & ($\alpha$)  &
$\left\{\begin{array}{l}
(A\drop \{v_*\})\cup\{hi,jk\} \\
(A\drop \{v_*\})\cup\{hj,ik\}
\end{array}\right.$
\\ \hline
(b1) & 2:2 & ($\beta$)\&($\gamma$1) &
$\left\{\begin{array}{l}
(A\drop \{v_*\})\cup\{hi,jk\} \\
(A\drop \{v_*\})\cup\{hj,ik\}
\end{array}\right.$
\\ \hline
(b2) & 1:1 & ($\gamma$2) & $(A\drop \{v_*\})\cup\{hi,jk\}$ \\ \hline
     & 0:1 & ($\gamma$3) & sign-reversing involution \\ \hline
\end{tabular}
\end{center}
\caption{Summary of Case 4(ii).}
\end{table}
The situations in which $(A\cap K,B\cap K)=(\{k\},\{h,i,j\})$ are handled analogously
after exchanging $A$ with $B$.\\

The next step is to identify which pairs $(A',B')$ contributing to equation (3) for
the $H$ graphs are produced by the above construction.  We refer to the edges
$ab$, $ac$, $ad$, $bc$, $bd$, and $cd$ as \emph{new edges}, so each $H$ graph has
two new edges.  The case analysis for these pairs $(A',B')$ is:\\
(i)\ the two new edges of $H$ are in different sets from the pair $(A',B')$;\\
(ii)\ the two new edges of $H$ are in the same set from the pair $(A',B')$.\\
We describe in detail the situation in case (ii) in which both new edges
are in $A'$.  (The other situation is analogous.)\\
(ii-$\alpha$)\ both $A'$ and $B'$ are spanning trees;\\
(ii-$\beta$)\ both $A'$ and $B'$ are forests with two components, and two of
the vertices $\{a_*,b_*,c_*,d_*\}$ are in the same component of $A'$ minus the new edges.\\
(ii-$\gamma$)\ both $A'$ and $B'$ are forests with two components, and each of
the vertices $\{a_*,b_*,c_*,d_*\}$ is in a different component of $A'$ minus the new edges.\\
Case (ii-$\gamma$) is further divided as follows:\\
(ii-$\gamma$1)\ neither of the new edges is in the cycle $C'$ of $A'\cup\{e,f\}$;\\
(ii-$\gamma$2)\ exactly one of the new edges is in the cycle $C'$ of $A'\cup\{e,f\}$;\\
(ii-$\gamma$3)\ both of the new edges are in the cycle $C'$ of $A'\cup\{e,f\}$.\\

The inverse of the above construction is described as follows.

In case (i), if the new edges of $H$ are $hi\in A'$ and $jk\in B'$ then put
$$A^\circ=(A'\drop \{hi\})\cup\{h,i\}\ \ \ \mathrm{and}\ \ \ B^\circ=(B'\drop \{jk\})\cup\{j,k\}.$$

In case (ii) we describe the situation in which the new edges of $H$ are $hi$ and $jk$,
both in $A'$.  The other situation is analogous upon exchanging $A'$ with $B'$.

In case (ii-$\alpha$), since $A'$ is a spanning tree there is a unique path in $A'$
with the new edges $hi$ and $jk$ as end-edges.  Let $i_*$ and $j_*$ be the end-vertices of this path,
and let $h_*$ and $k_*$ be the other ends of the new edges.  We associate two pairs with
$(A',B')$ in this case, by
$$A^\circ=(A'\drop \{hi,jk\})\cup\{h,i,j\}\ \ \ \mathrm{and}\ \ \ B^\circ=B'\cup \{k\}.$$
and
$$A^{\circ\circ}=(A'\drop \{hi,jk\})\cup\{i,j,k\}\ \ \ \mathrm{and}\ \ \ B^{\circ\circ}= B'\cup \{h\}.$$

In case (ii-$\beta$), there is also a unique path in $A'$ with the new edges $hi$ and $jk$ as
end-edges.  We associate two pairs with $(A',B')$ in this case by the same construction as in case
(ii-$\alpha$).

In case (ii-$\gamma$1), there are two paths in $A'\cup\{e,f\}$ with the new edges $hi$ and $jk$ as
end-edges.  Notice that these two paths have the same end-vertices.  Thus, we let $i_*$ and $j_*$ be
the common end vertices of these paths, we let $h_*$ and $k_*$ be the other ends of the new edges,
and we associate two pairs with $(A',B')$ in this case by the same construction as in case
(ii-$\alpha$).

In case (ii-$\gamma$2), let $hi$ be the new edge of $H$ in the cycle $C'$ of $A'\cup\{e,f\}$, and let
$jk$ be the new edge of $H$ not in $C'$.  Since we are in case (ii-$\gamma$), 
the edges $hi$ and $jk$ are in different components of $A'$.  Let $k_*$ be the end of $jk$ that
is in the connected component of $(A'\cup\{e,f\})\drop \{jk\}$ which contains $C'$.  We associate
one pair with $(A',B')$ in this case, by
$$A^\circ=(A'\drop \{hi,jk\})\cup\{h,i,j\}\ \ \ \mathrm{and}\ \ \ B^\circ=B'\cup\{k\}.$$

In each case except (ii-$\gamma$3) we associate one or two pairs $(A^\circ,B^\circ)$,
$(A^{\circ\circ},B^{\circ\circ})$ to each pair $(A',B')$ contributing to (3) for the three $H$ graphs.
In each case, pairs that are associated with one another are of the same type.  By considering the
constructions given for Case 4(ii) in detail, one sees that they give correspondences as in Table 2.  
Thus, the pairs contributing to (3) for $G$ are equinumerous (and of the same types) as the pairs
contributing to (3) for the $H$ graphs, except for the pairs in case (ii-$\gamma$3).  Comparing this with
equations (6) and (7), we see that $\LDIFF$ is the number of pairs in case (ii-$\gamma$3)
of types $\X^+\X^-$ or $\X^-\X^+$, and that $\RDIFF$ is the number of pairs in case
(ii-$\gamma$3) of types $\X^+\X^+$ or $\X^-\X^-$.

The following sign-reversing involution on the set of pairs in case (ii-$\gamma$3) shows that
$\LDIFF=\RDIFF$ and completes the analysis of Case 4, the induction step, and the proof.
Consider any pair $(A,B)$ in case (ii-$\gamma$3), with new edges of $H$ being $hi$ and $jk$,
both in $A$.  (The situation in which both new edges are in $B$ is analogous.)
Notice that since we are in case (ii-$\gamma$), the two new edges are in different components of
$A$.  Exactly one of the sets
$$
(A\drop\{hi,jk\})\cup\{hj,ik\}
\ \ \ \mathrm{or}\ \ \
(A\drop\{hi,jk\})\cup\{hk,ij\}
$$
is also a forest $F$ with two components such that both $F\cup \{e\}$ and $F\cup \{f\}$ are spanning trees.
Call this set $A^\natural$, and note that $(A^\natural)^\natural=A$.  The sign-reversing involution
is given by
$$(A,B)\ \longleftrightarrow (A^\natural,B)$$
as illustrated in Figure 7.
\pix{1.0}{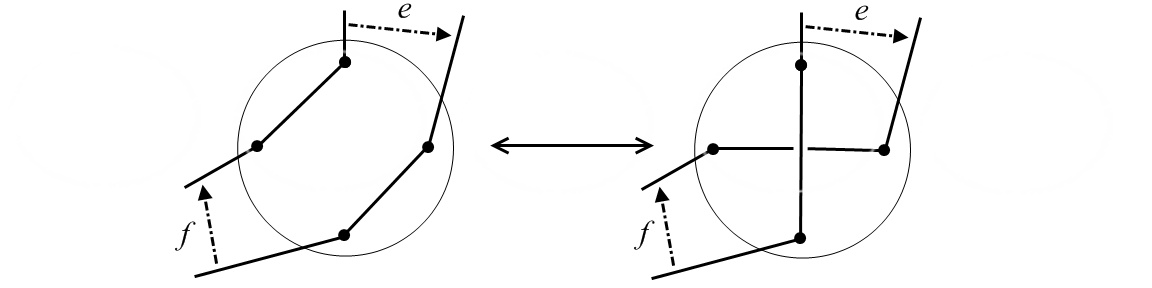}{The sign-reversing involution of Case 4(ii-$\gamma$3).}
Thus, in case (ii-$\gamma$3) the pairs of types  $\X^+\X^-$ or $\X^-\X^+$ are equinumerous with
the pairs of types $\X^+\X^+$ or $\X^-\X^-$.  That is, $\LDIFF=\RDIFF$.  This completes the proof.
\hfill $\Box$


\begin{thebibliography}{BSST}

\bibitem{BB} N. Balabanian and T.A. Bickart,
``Electrical Network Theory,'' Wiley, New York, 1969.

\bibitem{BSST} R.L. Brooks, C.A.B. Smith, A.H. Stone, and W.T. Tutte,
\emph{The dissection of rectangles into squares},
Duke Math. J. \textbf{7} (1940), 312-340.

\bibitem{Cha} S. Chaiken,
\emph{A combinatorial proof of the all minors matrix tree theorem},
SIAM J. Alg. Disc. Methods  \textbf{3} (1982), 319-329.

\bibitem{Ch1} Y.-B. Choe, ``Polynomials with the Half-Plane
Property and Rayleigh Monotonicity,'' Ph.D. Thesis, University of 
Waterloo, 2003.

\bibitem{Ch2} Y.-B. Choe, \emph{A combinatorial proof of the Rayleigh formula
for graphs}, Discrete Math., to appear.

\bibitem{Ch3} Y.-B. Choe, \emph{Sixth-root of unity matroids are Rayleigh},
in preparation.

\bibitem{COSW} Y.-B. Choe, J.G. Oxley, A.D. Sokal, and D.G. Wagner,
\emph{Homogeneous polynomials with the half-plane property},
Adv. Appl. Math. \textbf{32} (2004), 88-187.

\bibitem{CW} Y.-B. Choe and D.G. Wagner, \emph{Rayleigh matroids},
Combin. Probab. and Comput. \textbf{15} (2006), 765-781.

\bibitem{CH} J. Cibulka and J. Hladk\'y,
\emph{Elementary proof of Rayleigh formula for graphs},
in ``Proceedings of SVO\v{C} 2007.'' (7 pp), and
\texttt{http://arxiv.org/abs/0803.4395}.

\bibitem{FM} T. Feder and M. Mihail, \emph{Balanced matroids},
in ``Proceedings of the 24th Annual ACM (STOC)'', Victoria B.C.,
ACM Press, New York, 1992.

\bibitem{Ki} G. Kirchhoff, 
\emph{\"Uber die Aufl\"osung der Gleichungen, auf welche man bei der
Untersuchungen der linearen Vertheilung galvanischer Str\"ome gef\"uhrt wird},
Ann. Phys. Chem. \textbf{72} (1847), 497-508.

\bibitem{Max} J. Clerk Maxwell, ``Electricity and Magnetism,'' Clarendon Press,
Oxford, 1892.
 
\end{thebibliography}
\end{document}